\documentclass[11pt]{article}

\usepackage{amscd,amsmath,amsxtra,amssymb,latexsym,amsfonts,graphics}
\usepackage[usenames,dvipsnames]{color}

\usepackage[small]{caption}

\usepackage{ifthen}
\newboolean{PrivateVersion}
\setboolean{PrivateVersion}{false} 

\newcommand{\private}[1]{ \ifthenelse{\boolean{PrivateVersion}}
{#1}{} }

\newtheorem{Thm}{Theorem}[section]
\newtheorem{Lemma}[Thm]{Lemma}

\newtheorem{Assumption}[Thm]{Assumption}

\newtheorem{Definition}[Thm]{Definition}

\newcommand{\R}{{\mathbb R}}

\newcommand{\C}{{\mathbb C}}
\newcommand{\Z}{{\mathbb Z}}

\renewcommand{\Im}{{\rm Im \ }}
\renewcommand{\Re}{{\rm Re \ }}

\newcommand{\xu}{{\underline{x}}}
\newcommand{\pravo}{{\rm right}}
\newcommand{\levo}{{\rm left}}

\makeatletter
\newif\if@caption@empty
\newcommand{\captionempty}
{%
  \@caption@emptytrue
  \caption{}%
}
\renewcommand{\@makecaption}[2]
{%
  \centering
  \itshape
  \rm #1%
  \if@caption@empty
    \global\@caption@emptyfalse
  \else
    \textbf{:} #2%
  \fi
}
\makeatother

\newcommand{\Ln}{{\rm Ln \ }}

\renewcommand{\Im}{{\rm Im \ }}

\begin{document}

\title{ Shatalov-Sternin's construction of complex WKB solutions and the choice of integration paths. 
\ifthenelse{\boolean{PrivateVersion}}
{\textcolor{blue}{ \\ !!!PRIVATE VERSION!!!}}{}  }

\author{Alexander GETMANENKO \\
{\footnotesize Institute for the Physics and Mathematics of the Universe,} \\
{\footnotesize The University of Tokyo, 5-1-5 Kashiwanoha, Kashiwa, 277-8568, Japan} \\ 
{\tt \small Alexander.Getmanenko@ipmu.jp}}

\maketitle

\begin{abstract}
We re-examine Shatalov-Sternin's proof of existence of resurgent solutions of a linear ODE. In particular, we take a closer look at the ``Riemann surface" (actually, a two-dimensional complex manifold) whose existence, endless continuability and other properties are claimed by those authors. We present a detailed argument for a part of the ``Riemann surface"  relevant for the exact WKB method. \\ The present text is the author's article {\tt arXiv:0907.2934} rewritten from a different perspective.
\end{abstract}

\section{Introduction.}

\subsubsection*{Resurgent analysis.}

Resurgent analysis is a method of studying hyperasymptotic expansions 
\begin{equation}  \sum_{k,j} e^{-c_k/h} a_{k,j}h^j , \ \ \ h\to 0+  \label{hyperasexpn}  \end{equation}
and those of similar kind by treating such expansions as asymptotics obtained from a Laplace integral 
\begin{equation} \int_\gamma \Phi (s) e^{-s/h} ds, \label{LaplaceTransf} \end{equation} 
where $\Phi$ is a ramified analytic function in the complex domain with a discrete set of singularities and $\gamma$ is an infinite path on the Riemann surface of $\Phi$. The crucial observation is that the terms of (\ref{hyperasexpn}) can be recovered from studying the singularities of $\Phi$, see ~\cite{V83}, ~\cite{E81}, ~\cite{CNP}, ~\cite{DP99}, as well as ~\cite{G} for this author's preferred terminology. 

The methods of resurgent analysis have been used, in particular, to study asymptotics of solutions of linear ODE with a small parameter, especially the Schr\"odinger equation in the semiclassical approximation, see, e.g. ~\cite{DDP97}; this technique is a refinement of what is known as the {\it complex WKB method}. More specifically, consider an equation of the type
\begin{equation} -h^2 \partial_x^2 \varphi(h,x) + V(x) \varphi(h,x) = 0 \label{SchroeEq} \end{equation}
where $x$ ranges over $\C$, $h$ is a small complex asymptotic parameter, and $V(x)$ is an entire function often assumed to be a polynomial.  Under the transformation (\ref{LaplaceTransf}), this equation becomes an equation on an unknown ramified analytic function of two variable $\Phi(s,x)$ of the form
\begin{equation} - \partial_x^2 \Phi(s,x) + \partial_s^2 V(x)\Phi(s,x) \ = \ 0. \label{MainEqu} \end{equation}
The equation \eqref{MainEqu} only needs  to be satisfied modulo functions that are entire with respect to $s$ for every value of $x$ since such functions correspond to zero under a properly (~\cite[Pr\'e I.2]{CNP}) understood Laplace transform \eqref{LaplaceTransf}.  Since the beginnings of resurgent analysis in the early 1980s there has been no real doubt that (\ref{MainEqu}) possesses two linearly independent (in an appropriate sense) solutions that are endlessly analytically continuable with respect to $s$ and satisfy certain growth conditions at infinity.

The manifold on which $\Phi(s,x)$ is defined is usually quite complicated. In the special cases when $V(x)=x$ and $V(x)=x^2$, the function $\Phi(s,x)$ can be written down by an explicit formula and $\varphi(h,x)$ is expressible in terms of Airy or Weber function, see \cite{J94}. For more complicated potentials, say, when $V(x)$ is a generic polynomial of degree $\ge 4$, the function $\Phi(s,x)$ is expected to be defined on a highly transcendental manifold, see \cite{DDP93} and \cite{D92}: if for a fixed $x$ one projects {\it all} singularities on {\it all} sheets of the Riemann surface of $\Phi(s,x)$ to the complex plane of $s$, one expects to obtain an everywhere dense set. Thus, there is no hope that the manifold in question is a universal cover of $\C^2$ minus a discrete family of complex curves.

\begin{figure} \includegraphics{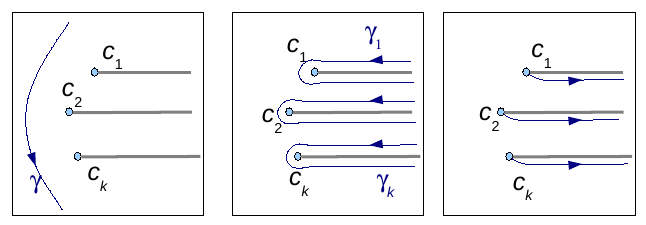} \caption{Deformation of the integration contour and the calculation of the hyperasymptotic expansion of $\varphi(h,x)$} \label{Paper3p3}
\end{figure}

Singularities of $\Phi$ and the precise structure of the manifold on which $\Phi$ is defined are important because they allow us to obtain the hyperasymptotic expansion of $\varphi(h,x)$ for $h\to 0+$ as follows (cf. ~\cite[p.218]{V83}, \cite{CNP}).  Fix $x$ and identify one of the sheets of the Riemann surface of $\Phi(x,s)$ with a complex plane of $s$ minus countably many cuts $c_1+\R_{\ge 0}, c_2+\R_{\ge 0}, ..., c_k+\R_{\ge 0}$ in the positive real direction.    Draw an infinite integration path $\gamma$ in $\C$ to the left of $c_1,c_2,...,c_k,...$, fig.\ref{Paper3p3},left,  so that, at least morally, $\varphi(h,x)=\int_\gamma \Phi(s,x) e^{-s/h} ds$. Using analyticity of $\Phi(s,x)$ and under appropriate conditions on its growth at infinity one can push the integration contour $\gamma$ to the right and rewrite 
$$ \varphi(h,x) \ = \ \sum_k \int_{\gamma_k} \Phi(s,x) e^{-s/h} ds, $$
where infinite integration paths $\gamma_k$ ``hang" on the singularities $c_k$, fig.\ref{Paper3p3},middle. Finally, one deforms each $\gamma_k$ so that both infinite branches lie on different sheets of the Riemann surface right on top of each other, and rewrites 
\begin{equation} \int_{\gamma_k} \Phi(s,x) e^{-s/h} ds \  = \ \int_{[c_k, c_k+\infty)}  (\Delta_{c_k}\Phi(s,x) )e^{-s/h} ds, \label{eq442} \end{equation}
where $\Delta_{c_k}\Phi$ denotes the jump of $\Phi$ across the cut starting at $c_k$. The integrals on the R.H.S. of \eqref{eq442} are taken along semi-infinite real analytic paths similar to those on fig.\ref{Paper3p3},right. The asymptotic expansions of these integral can now be calculated using Watson's lemma and combined to a hyperasymptotic expansion \eqref{hyperasexpn}. 

In ~\cite{CNP}, ~\cite{ShSt}, the following point of view is developed. For each fixed $x$, $\Phi(s,x)$ as a function of $s$ is assumed in the beginning to be a holomorphic function on a sectorial neighborhood of infinity $\Omega_0=\{ s\in \C \ : \ \arg s \in (\frac{\pi}{2}-\beta, \frac{3\pi}{2}+\beta;\  |s|>N\}$ for some $\beta>0$ and $N>0$; the contour $\gamma$ appearing in \eqref{LaplaceTransf} is a contour along the boundary of $\Omega_0$. It is then assumed that for a discrete subset $\{c_1,..,c_k,.. \}\subset \C\backslash \Omega_0$, the function $\Phi(s,x)$ has an analytic continuation to the set $\Omega=\C\backslash \bigcup_k (c_k+\R_{\ge 0})$; this $\Omega$ is called {\it the first sheet} of the Riemann surface of $\Phi(s,x)$, and the points $ c_k$, $k=1,2,...$, are called the {\it the first sheet singularities} of $\Phi$. The Riemann surface of $\Phi(s,x)$ for every fixed $x$ is the Riemann surface of the analytic continuation of $\Phi(s,x)$ as a holomorphic function on $\Omega_0$. It is important that in order to  obtain $\Delta_{c_k} \Phi(s,x)$ in \eqref{eq442} as an analytic function of $s$, we define it as  $\Phi(s',x)-\Phi(s'',x)$ where $s',s''$ belong to the different sheets of the Riemann surface of $\Phi(x,s)$ and project to the same point of $s\in\C$; we need therefore an analytic continuation of $\Phi$ beyond the first sheet at least near the cuts $c_k+\R_{\ge 0}$. 
 
While the position of the singularities of $\Phi(s,x)$ is important for the calculation of the asymptotics, there is a good intuition where these singularities are located. Given an initial point $x_0$ and  two ramified analytic functions $f_0(s), f_1(s)$, let $\Phi(s,x)$ solve the Cauchy problem $\Phi(s,x_0)=f_0(s)$, $\frac{\partial}{\partial x}\Phi(s,x_0)=f_1(s)$ for the equation \eqref{MainEqu}. The general philosophy of PDE suggests that the singularities of the initial conditions should propagate along the integral curves of the vector fields $\frac{\partial}{\partial x} \pm \sqrt{V(x)}\frac{\partial}{\partial s}$. Using this intuition, Voros ~\cite{V83} studied the Stokes phenomenon -- appearance and disappearance of singularities from the first sheet of $\Phi(s,x)$ as $x$ varies, and described its consequences (``connection formulas") for the hyperasymptotic expansions of $\varphi(h,x)$.

Since so much relies on the properties of singularities and analytic continuation of $\Phi(s,x)$, proving that \eqref{MainEqu} has an endlessly analytically continuable solution is an important foundational question.  
The present work is a step in this direction.

\subsubsection*{Literature review}

The literature on this subject is extremely vast, so we can hope to at most indicate {\it some} sources which reflect the state of the field and main developments.

The problem of existence and singularities of complex-analytic solutions $\Phi$ of \eqref{MainEqu} appear in numerous classical works, notably ~\cite{Le}, ~\cite{Ha} and their sequels, but the solutions are shown to exist only locally, and the results do not guarantee existence of the analytic continuation of $\Phi(s,x)$ to the values of $x$ far away from an initial point $x_0$ where the Cauchy data are given. 

From ~\cite{DP99} we learned about the existence of a preprint ~\cite{E84} containing a sketch of a construction of endlessly continuable solutions $\Phi$ satisfying \eqref{MainEqu}, but at least according to ~\cite{DP99}, not all details are clear in that sketch. 

Lacking a general statement, one could still work out examples of potentials $V$ for which the function $\Phi$ can be given by a more or less explicit formula and singularities of $\Phi$ are possible to analyze from that explicit representation, see e.g. the easiest examples in  ~\cite{J94} and much more complicated one in the recent article ~\cite{FS}.  

In the terminology of resurgent analysis, the function $\Phi$ appearing in \eqref{MainEqu} is the ``major" of $\varphi$ appearing in \eqref{SchroeEq}. Many authors prefer to take a somewhat different Laplace integral and work with ``minors"; there is a technology of translating statements between the two setups, ~\cite{CNP}.  Working with minors, the authors of ~\cite{DLS93} present a proof that we expect to imply the existence of $\Phi(s,x)$ for values of $s$ on the first sheet minus the cuts and for $x$ confined to a region where no Stokes phenomenon occurs. 

The monograph ~\cite[Ch.3.1]{ShSt} and numerous works by the same authors, e.g. ~\cite{SS93}, ~\cite{SSS97},  contain another approach to the proof of existence of endlessly continuable solutions of \eqref{MainEqu} and of similar equations of higher order. From the parts of the argument that we were able to understand, 
the approach seems very natural and attractive. Discussion  of ~\cite{ShSt}'s proof is the content of this article. 

The topic has remained in the focus of many researchers. It may have been one of the motivations for development of the mould calculus, cf. ~\cite{Sa} and references therein.

Meanwhile the Kyoto school has been working on the idea of transforming the
Schr\"odinger equation with an arbitrary potential $V(x)$ to appropriately chosen canonical models, e.g. Airy,  Weber, Whittaker equations, e.g. ~\cite{AKT91},  ~\cite{KKKT10}; the language of ``minors" is used by these authors.
A breakthrough was announced in the autumn of 2010 by Kamimoto and Koike. Their result is expected to describe the first sheet singularities of $\Phi(s,x)$ as a function of $s$, as long as $x$ is close to a simple zero of a very general potential $V(x)$.  

Not only \eqref{SchroeEq}, but also other similar equations have been studied by means of complex WKB method; respectively, different equations in the Laplace-transformed picture take the place of \eqref{MainEqu}. E.g., higher order ODEs were studied semi-heuristically in ~\cite{AKSST05}, ~\cite{H08}, or rigorously in ~\cite{NNN91}; the  first order difference equations with a small parameter were studied in ~\cite{CG08}. 

In the present article we are re-examining certain details of the Shatalov-Sternin's proof. The idea of the argument presented in ~\cite{ShSt} differs significantly from what the approach of the Kyoto school and from that of other authors. Even in  view of the results announced by Kamimoto and Koike it remains important, for our understanding of the subject as well as for possible extensions and generalizations, to clarify the status of ~\cite{ShSt}'s very natural-looking argument. 

At the time when this version of the article is written, its ideas have been already used in ~\cite{GT}. 

\subsubsection*{Contribution of this article.}

In ~\cite{ShSt}, Sternin and Shatalov solve \eqref{MainEqu} by reducing it to an integral equation and obtaining a resolvent. This method is classical in the theory of differential equations; an anonymous referee pointed out that it was already used, in a slightly different setup, in ~\cite[Ch.2]{BB74}. In other words, the authors of ~\cite{ShSt} represent a solution $\Phi(s,x)$ in terms of an infinite series 
\begin{equation} \Phi(s,x) = \sum_{n=0}^\infty \Psi_n(s,x) \label{Jan6} \end{equation}
where $\Psi_n(s,x)$ is, morally, the result of an $n$-fold application of some integro-differential operator to a ``0-th order approximation" $f(s)$. 
The actual formulas will be recalled in section \ref{SSconSec}.

Having formally obtained an expression \eqref{Jan6}, Sternin and Shatalov set out to prove that a) all functions $\Psi_n(s,x)$ are defined on the same endlessly continuable manifold of complex dimension two (which is still called a ``Riemann surface"), and that b) the series converges on compact sets of this ``Riemann surface".

In ~\cite[Prop.3.1, pp.204-207]{ShSt}, the construction of the ``Riemann surface" takes only three pages and is presented very intuitively; however, once we wanted to make a precise sense of how exactly the ``Riemann surface" is described and how exactly all functions $\Psi_n$ can be analytically continued to it by which specific deformations of integration contours, we found ourselves dealing with a rather complex situation.

In this article,  we restrict ourselves to constructing an open piece ${\cal S}$ of the ``Riemann surface". As a bit of an oversimplification, let us say that over each point $x$ in an appropriate region of the complex plane, the fiber of ${\cal S}$ is a complex plane minus finite number of rays, ``cuts",  in the positive real direction.

Thus, the statement and the proof of the following theorem are intended to make precise some things which we could not find in ~\cite{ShSt}. 

\begin{Thm} \label{ContribTh} For $V(x)$ satisfying assumptions of section \ref{WhoIsVx}, the countably many functions  \eqref{vNsimpl} possess an analytic continuation to the 2-dimensional complex manifold ${\cal S}$ defined in section \ref{FibOfS}. \end{Thm} 

A word of caution: The functions $\Psi_n(s,x)$ appearing in \eqref{Jan6} are more complicated than functions \eqref{vNsimpl}, but it will be obvious that the theorem implies that $\Psi_n$ also analytically continue to ${\cal S}$.

Here is what remains outside the scope of this article. The series \eqref{Jan6} is very likely to converge uniformly on compact subsets of ${\cal S}$. Unfortunately, in ~\cite[(3.14)]{ShSt} the derivative in the integrand of ~\eqref{Rj} is missing, and those authors end up proving convergence of a wrong and much better behaving series. A more delicate study of convergence will need to be performed in the future. The current paper makes the question more well-defined: before we study convergence of the series \eqref{Jan6} at a point $(s,x)$ of ${\cal S}$, we need to know first {\it how exactly} the functions $\Psi_n$ are analytically continued to the point $(s,x)$.
If the convergence is shown, that will provide an alternative both to the approach announced by Kamimoto and Koike and to the method of ~\cite[GT]. 

Let us briefly indicate what is involved in the proof of theorem \ref{ContribTh}. As the ramified analytic functions  \eqref{vNsimpl} of variables $(s,x)$ are iterations of two integro-differential operators $R_1$ and $R_2$, in order to analytically continue these functions to a point $(s,x)$ we need to appropriately define two integration paths (one for $R_1$ and one for $R_2$) leading from $(s_0,x_0)$ to $(s,x)$; here $x_0$ is some fixed initial point and $s_0$ depends on $s$ and $x$. First we treat the case when $x$ is in the same Stokes region as $x_0$, and then describe, starting from section \ref{modelcases}, a method that allows us to draw the integration paths for $x$ belonging to further and further Stokes regions. As we take $x$ in Stokes regions further and further away from $x_0$, there appear more and more obstacles to drawing an integration path from $(s_0,x_0)$ to $(s,x)$; points $(s,x)$ that cannot be reached by an integration path give rise exactly to
 the  singularities of ${\cal S}$
 predicted by Voros.

\paragraph*{The related paper {\tt arXiv:0907.2934}.} The open piece of the ``Riemann surface" ${\cal S}$ which we are constructing in the present paper is insufficient for the deformations of the integration contour that we need in \eqref{eq442}. 

In {\tt arXiv:0907.2934}, we are constructing, in a special situation, a larger ${\cal S}$ so that its fibers ${\cal S}_x$ are comprised of a the first sheet (i.e. the complex plane with finitely many cuts)  and small ``flaps"  attached on the sides along each cut; if we knew how to prove the convergence of \eqref{Jan6}, that would suffice to fully justify the procedure of \eqref{eq442}. Constructing this larger ${\cal S}$ is done similarly to the present paper, but requires much heavier notation and leads to a less crisp result.

\section{Shatalov-Sternin's construction.} \label{SSconSec}

The purpose of this section is to review the content of  ~\cite[pp.198-204]{ShSt} in the special case of the one-dimensional Schr\"odinger equation
\begin{equation}  [-h^2 \partial^2_x + V(x)]\varphi(h,x) = 0, \label{SchroeEq1} \end{equation}
where the variable $x$ takes values in $\C$ and $V(x)$ is an entire function. 

To describe the Laplace-transformed version of \eqref{SchroeEq1}, consider the following operation on the equivalence classes of germs of analytic functions at a point $(s_0,x_0)\in \C^2$ modulo functions entire with respect to $s$ for every $x$: 
$${\hat h}\Phi(s,x) \ := \  \partial_s^{-1} \Phi(s,x) = \int^s_{s_*(x)} \Phi(s',x)ds',$$
where the starting point of the integration $s_*(x)$ may depend on $x$. Changing $s_*(x)$ will change the result by a function depending only on $x$. 

In this notation, the Laplace transform  \eqref{LaplaceTransf} turns \eqref{SchroeEq1} into 
\begin{equation} -{\hat h}^2 \partial_x^2 \Phi(s,x) + V(x)\Phi(s,x) \ = \ 0 \label{MainEqu1} \end{equation}
which has to be satisfied modulo functions that are entire with respect to $s$ for every $x$. We would like to find solutions $\Phi$ of \eqref{MainEqu1} that are holomorphic functions on a complex two-dimensional manifold ${\cal S}$ endowed with a locally biholomorphic projection $\Pi$ to $\C^2$ with coordinates $(s,x)$. We would also like, for every $\xu\in\C$, the connected components of $\Pi^{-1}(\{ (s,\xu):s\in \C\})$ to be endlessly continuable Riemann surfaces in the sense of resurgent analysis, e.g., ~\cite[R\'es I]{CNP}. In fact, ~\cite{ShSt} use the concept of a ``ramified analytic function" of several complex variables; we will replace it by a clearer notion of ``a germ of an analytic function" except in philosophical statements.

The Cauchy-Kowalewskaya theorem, e.g.~\cite[Th.3.1.1]{Sch}, or the related results of ~\cite{Le} and ~\cite{Ha}, for this equation fall far short of the statement that we need. Indeed, for the equation \eqref{MainEqu1} with an initial condition, say, $\Phi(s,x_0)=\frac{1}{2\pi i s}$, $\frac{\partial}{\partial x}\Phi(s, x_0)=0$ (corresponding to $\varphi(h,x_0)=1$, $\frac{\partial}{\partial x}\varphi(h,x_0)=0$) one would only get existence of solution $\Phi(s,x)$ in a small polydisc centered at $(s_0,x_0)$ for $s_0\ne 0$, and the size of that polydisc is hard to increase. Therefore a more explicit construction of $\Phi$ is proposed. 

Fix a point $x_0$ such that $V(x_0)\ne 0$ and a determination $p(x)$ of $\sqrt{V(x)}$ in a neighborhood of $x_0$. Let $p_1(x)=-p_2(x)=p(x)$; let further  $S_j(x) = \int_{x_0}^x p_j(y) dy$, $j=1,2$, and $S(x)=S_1(x)$. 
In this notation, the operator $-\hat h^2 \partial_x^2 +V(x)$ on the L.H.S. of \eqref{MainEqu1} can be rewritten as 
\begin{equation} \left(p^2(x) [-\frac{1}{p(x)} \hat h\partial_x - 1] - {\hat h}p'(x) \right)\left(\frac{1}{p(x)}{\hat h}\partial_x - 1\right) - {\hat h}p'(x). \label{splitOrd1} \end{equation} 

We will be able to make use of this representation once we are able to invert the operators $\pm \frac{1}{p(x)} \hat h \partial_x -1$. Namely, consider an equation
\begin{equation} [\frac{1}{p_j(x)} \hat h \partial_x -1] u(s,x) = b(s,x), \label{e308} \end{equation}
as an equation of germs at $(s_0,x_0)\in \C^2$ of analytic functions of $(s,x)$  modulo functions depending only on $x$. Then \eqref{e308} is satisfied by 
$$ u(s,x) \ = \ R_j b(s,x) + f(s+S_j(x)), $$
where $f(s)$ is any germ of an analytic function near $s_0$ and the operator $R_j$  is defined by the formula
\begin{equation} (R_j G)(s,x) = \int_{x_0}^x (D_1G)(s+S_j(x)-S_j(y),y) p_j(y) dy.  \label{Rj} \end{equation}
Here $D_1$ stands for the derivative of the function with respect to the first argument. We consider $R_j$ as acting on germs of analytic functions $G(s,x)$ at a point $(s_0,x_0)$. In ~\cite{ShSt} this derivative is missing.

Let us start looking for a solution \eqref{MainEqu1} in the form
\begin{equation} \Phi(s,x) = R_1 \Phi_1(s,x) + f_1(s+S_1(x)).  \label{e327} \end{equation}
Substituting \eqref{e327} into \eqref{MainEqu1} and using the expression \eqref{splitOrd1}, we have
\begin{equation} \left\{\left(p^2(x) [-\frac{1}{p(x)} {\hat h}\partial_x - 1] - {\hat h}p'(x) \right) - {\hat h}p'(x) R_1\right\} \Phi_1 \ = \
 - {\hat h} p'(x) f_1(s+S_1(x)).  \label{e330} \end{equation}
Looking for a solution of \eqref{e330} in the form
$$ \Phi_1(s,x) = R_2\Phi_2(s,x) + f_2(s+S_2(x)), $$
we obtain
{ \small
$$ \left[ 1 - {\hat h}\frac{p'(x)}{p^2(x)}\{ R_2  + R_1 R_2 \} \right] \Phi_2  \ = \  - {\hat h} \frac{p'(x)}{p^2(x)}  \{ (1+R_1)  f_2(s+S_2(x)) + f_1(s+S_1(x)) \}.  $$
}

Formally, the last equation has a solution
\begin{equation} \Phi_2(s,x) \ = \  \sum_{j=0}^{\infty} (-1)^j {\hat h}^{j+1} [ (-\frac{p'(x)}{p^2(x)})(1+R_1) R_2 ]^j g_0(s,x) , 
\label{vNseries} \end{equation} 
where 
$$ g_0(s,x) \ = \ - {\hat h} \frac{p'(x)}{p^2(x)}  \{ (1+R_1)  f_2(s+S_2(x)) + f_1(s+S_1(x)) \}. $$
 On the R.H.S. of \eqref{vNseries} we see an infinite series of germs of analytic functions; only its partial sums are mathematically well-defined at this stage.

Assume that we are able to prove that the series on the right hand side of \eqref{vNseries} converges both for the choice a) $f_1(s)=\Ln s$, $f_2=0$, and for the choice b) $f_1=0$, $f_2(s)=\Ln s$, and in both cases defines analytic functions $\Phi_2(s,x)$ and $\Phi(s,x)$ on a sufficiently large complex two-dimensional manifold. Then we can perform a Laplace integral as in \eqref{eq442}; as a result, we expect to obtain two formal WKB solutions of \eqref{SchroeEq1} for $x$ in a neighborhood of $x_0$, namely $A_+(h,x) e^{S(x)/h} + A_-(h,x)e^{-S(x)/h}$ for the choice a), and $B_+(h,x) e^{S(x)/h} + B_-(h,x)e^{-S(x)/h}$ for the choice b). Here $A_\pm(h,x), B_\pm(h,x)$ are expected to be formal (actually, Gevrey) power series in $h$ with $x$-dependent coefficients. We expect further that two vectors $[A_+(h,x_0), A_-(h,x_0)]$ and $[B_+(h,x_0), B_-(h,x_0)]$ in $\C[[h]]^2$ will be linearly independent over $\C[[h]]$, thus yielding two linearly independent resurgent solutions of \eqref{SchroeEq1} in every reasonable definition of this notion.

The first task is therefore to construct a ``Riemann surface'' -- a two dimensional complex manifold on which all summands in the R.H.S. of \eqref{vNseries} are defined for the choices a) and b) from the previous paragraph. It is easy to see that an equivalent question is to construct a ``Riemann surface'' on which all functions 
\begin{equation} R_{j_k}...R_{j_2}R_{j_1} f(s,x) , \ \ \ j_i=1,2, \ \ \ k\ge 0  \label{vNsimpl} \end{equation}
are defined for $f(s,x)=\Ln(s\pm S(x))$. 

This is the question we are dealing with in this article. The second task would be to show that the infinite series converges on this ``Riemann surface". Unfortunately, a derivative in the integrand is missing in ~\cite{ShSt}'s definition of operators $R_j$ and we cannot suggest an easy way to repair their convergence argument, but hope to give (or read!) an alternative proof elsewhere.


\section{Analytic continuation and integration paths} \label{ACIP}

In section \ref{Notation} we are going to precisely define the ``Riemann surface" ${\cal S}$ to which we will then be able to analytically continue the functions \eqref{vNsimpl}.

Let $\C_s$, $\C_x$ denote the complex planes of the variables $s$, $x$, respectively.

 This section \ref{ACIP} exposes the main idea of this article; its content will make precise sense after reading section  \ref{Notation}.  For now we will think of ${\cal S}$ as some complex two-dimensional manifold with a locally biholomorphic projection $\Pi:{\cal S}\to \C_{s}\times\tilde {\cal O}$, where $\tilde {\cal O}$ is a complex one-dimensional manifold with a locally biholomorphic projection to $\C_x$.  We will freely use $(s,x)$ as local coordinates on ${\cal S}$.

\subsection{Reduction of the problem to construction of the integration paths.} \label{RedToPaths}

Recall that we denote  $p_1(x)=-p_2(x)=p(x)$,  $S_j(x) = \int_{x_0}^x p_j(y) dy$, $j=1,2$, and the operators $R_j$ were defined by \eqref{Rj} as operators acting on germs of analytic functions. 

It will be obvious from the construction of ${\cal S}$ that the functions $\Ln (s\pm S(x))$ have analytic continuations to ${\mathcal S}$. Existence of analytic continuation of all terms of \eqref{vNsimpl} to ${\mathcal S}$ will follow by induction from the following

\begin{Thm}  If $G(s,x)$ is defined and analytic on ${\mathcal S}$, then $R_jG$, $j=1,2$ have analytic continuations to ${\mathcal S}$. \label{RjGmainTh}
\end{Thm}


A detailed proof of this theorem will be given in section \ref{PfMresSec}. In this section \ref{RedToPaths} we will introduce some of the terminology used in the proof; in section \ref{AppearsStokes} we will informally  explain the idea on which the proof is based.

If $G(s,x)$ were an analytic function on the whole $\C\times \tilde {\cal O}$, we could define $(R_jG)(s,\xu)$ by the formula 
\begin{equation} (R_j G)(s,\xu) = \int_{x_0}^\xu (D_1G)(s+S_j(\xu)-S_j(y),y) p_j(y) dy\label{Rj1} \end{equation}
 where the integral is taken along {\it any} path from $x_0$ to $\xu$ in $\tilde{\mathcal O}$. Since, however, $G(s,\xu)$ is defined on a complicated manifold ${\cal S}$, we need to find for each $(s,\xu)\in{\cal S}$ a path $y(t)$ in $\tilde{\mathcal O}$ from $x_0$ to $\xu$ satisfying the following 

{\bf Definition.} We say that a path $y(t)$ in $\tilde {\cal O}$  {\it can be  lifted to ${\cal S}$ parallel to $-S_j$ with endpoint $(s,\xu)$} if  $(s+S_j(\xu)-S_j(y(t)),y(t))$ defines a path in ${\mathcal S}$.

Intuitively, this condition means that the point $(s+S_j(\xu)-S_j(y(t)),y(t))$  does not ``leave" ${\cal S}$ and does not hit any of its singularities.

We will call such a $y(t)$ {\it an integration path for $(s,\xu)$ and $R_j$}; let us stress that the choice of a path $y(t)$ depends on the point $s$ in the fiber ${\cal S}_\xu$ of ${\cal S}$ over $\xu$. If the integration paths $y(t)$ continuously depend on $(s,\xu)$, using them in \eqref{Rj1} yields an analytic function $R_jG(s,x)$; construction of  $R_jG$ from $G$ is thus reduced to finding a family of  integration paths continuously depending on $(s,\xu)$. \label{LiftingDefined}

\private{\textcolor{blue}{It seems that this paragraph is no longer needed.} Note that if $\xu$ is contained in a compact subset $K\subset\tilde{\mathcal O}$ and $x_0\in K$, then there is a constant $N_K\in \R$ such that any integration path $y(t)$ from $x_0$ to $\xu$ contained in $K$ will satisfy the desired property \textcolor{blue}{which property?} for any $s$ with $\Re s<N_K$. \textcolor{blue}{HEAVY!} }


In Section \ref{FibOfS} we will describe the fibers  ${\cal S}_x$ of ${\cal S}$ over every $x\in \tilde {\cal O}$; we will define, for each $x\in \tilde{\cal O}$, a list of {\it singularities} in ${\cal S}$ each of which will be of the form $s=S_j(x)+c$, $j=1,2$, $c\in \C$, for appropriate constants $c$.

For $U\subset{\mathcal S}_\xu$ \label{Ucomes} let us try to construct an integration path $y(t)$ which for any endpoint $(s,\xu)$, $s\in U$, can be lifted to ${\cal S}$ parallel to $-S_j$. Suppose $S_j(x)+c$, $c\in \C$, is one of the singularities of ${\cal S}_x$. We want to make sure that $s+S_j(\xu)-S_j(y(t))$ avoids the singularity $S_j(y(t))+c$, i.e. we want the equality 
$$ 2S_j(y(t)) \ = \ s+S_j(\xu)-c $$
to hold for no point $y(t)$ along the integration path and for no point $s\in U$. That is to say, we want the integration path $y(t)$ to avoid the set 
\begin{equation} S_j^{-1}\left( \frac{U+S_j(\xu)-c}{2} \right)\subset \tilde{\mathcal O}. \label{SjInverse} \end{equation}
 We will need to carefully keep track of the appropriate branches of the functions involved in this expression.


It will turn out a posteriori that the condition that $(s+S_j(\xu)-S_j(y(t)),y(t))$ does not coincide with any of the singularities of ${\cal S}_{y(t)}$  is enough to guide us through the choice of the integration paths $y(t)$ for the point $(s,\xu)$. Once a choice of an integration path $y(t)$ is proposed, it is an extra logical step to check that its lifting parallel to $-S_j$ stays within ${\cal S}$; this however will always be obvious by inspection and not mentioned explicitly.

When constructing integration paths $y(t)$ for $R_j$, we found it convenient to construct the parallel transport of the set $U\in {\mathcal S}_{\xu}$ by defining $U(y(t))=U+S_j(\xu)-S_j(y(t))$. Then, as $t$ varies, the set $U(y(t))$ and the singularities  of type $-S_j(y(t))+const$ move with respect to the $s$-coordinate parallel to each other, and differently from the singularities of type $S_j(y)+const$.  For this reason, we will call singularities of type $-S_j(y)+const$ {\it stationary} and the singularities of type $S_j(y)+const$ {\it moving} singularities. When the index $j$ changes, the roles of moving and stationary singularities reverse. 

\begin{Definition} \label{cbpted} We will say that an open set $U\in {\cal S}_\xu$ {\it can be parallel transported parallel to $-S_j$ within ${\cal S}$ along the path $y(t)$} if $U(y(t))=U+S_j(\xu)-S_j(y(t))\subset {\cal S}_{y(t)}$ for every point $y(t)$ on the path. \end{Definition}



\subsection{Appearance of the Stokes curves in the construction of the integration paths.} \label{AppearsStokes}

This subsection \ref{AppearsStokes} is written informally and included for illustrative purposes only; the precise argument in the rest of the paper does not logically depend on it.

As the functions $S_j$ enter into the definition of the operators $R_j$, it is natural to choose $\tilde {\cal O}$ from the introductory paragraph of section \ref{ACIP} in such a way that $\tilde {\cal O}\to \C_x$ factors through the universal cover of $\C\backslash V^{-1}(0)$ with the base point $x_0$.

Let us discuss the construction of $R_1G_0$ and $R_2G_0$ for the function $G_0(s,x)=\Ln(s+S(x))$ (compare to \eqref{vNsimpl}). This function $G_0(s,x)$ is naturally defined on a ``Riemann surface" ${\cal S}_0$ whose fiber over any $x\in \tilde{\cal O}$ is a universal cover of $\C_s\backslash \{ -S(x)\}$. All iterations $R^n_1G_0$, $n\ge 1$, are also defined on ${\cal S}_0$: arbitrary paths in $\tilde{\cal O}$ can be chosen as integration paths for defining $R_1G_0$ {\it for this specific function $G_0$}. On the contrary, the ``Riemann surface" ${\cal S}_2$ of $R_2G_0$ necessarily has (at least) an additional singularity at $s=S(x)$: for any integration path $y(t)$ from $x_0$ to $\xu$ for $R_2$ and $(s=S(\xu),\xu)$, the integrand of \eqref{Rj1} is singular for $y=x_0$ because $S(x_0)=0$ and  $G_0$ has a singularity at $(x_0,0)$. 

Thus, the common ``Riemann surface" ${\cal S}$ of all the functions \eqref{vNsimpl} necessarily has singularities at $s=S(x)$ and $s=-S(x)$ on its first sheet.

\begin{figure} \includegraphics{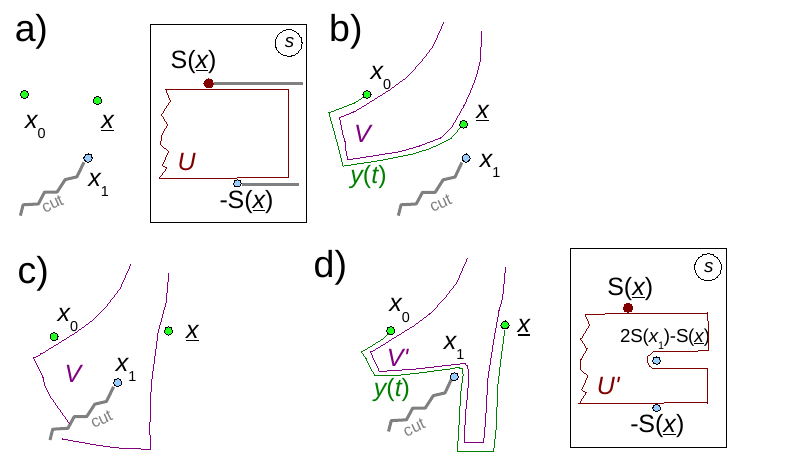}  \caption{Section \ref{AppearsStokes}.  
a) Projections of $x_0$ and $\xu$ to the complex plane of $x$; the set $U$ in the fiber of ${\cal S}_2$ over $\xu$; b) The set $V_\C$ and the integration path $y(t)$ in the case $\Im S(\xu) < \Im S(x_1)$; c) The set $V_\C$ in the case $\Im S(\xu) > \Im S(x_1)$; d) In the situation of c), the set $V'$ and the integration path $y(t)$ in the complex plane of $x$;  the set $U'$ in the fiber of ${\cal S}_2$ over $\xu$.
The branch cut starting from $x_1$ reminds us that the function $S(x)$ has a branch point at $x_1$.} 
\label{Paper3p44n} \end{figure}

 Suppose $x_1\in \C$ is a zero of $V(x)$, and all other zeros of $V(x)$ are far enough from $x_0$ so as not to affect our reasoning here; let $\Im S(x_1)>0$. Take a point $\xu\in {\cal O}$ such that $\Im S(\xu)>0$. Assume that the function $G_2(s,x)=(R_2G_0)(s,x)$ is defined on the Riemann surface ${\cal S}_2$ with singularities at $s=\pm S(x)$  on the first sheet. Let us study whether $R_1G_2(s,x)$ can be analytically continued to the set $U=\{ s\in \C : \Im [-S(\xu)]<\Im s <\Im S(\xu), \ \Re s < N \}$ identified with a subset in the fiber ${\cal S}_{2}$ over $\xu$, where $N\in \R$ is a large positive number, fig. \ref{Paper3p44n},a). 

The reasoning of \eqref{SjInverse} with $S(x)$ playing the role of $S_j(x)+c$ leads us to considering the set $V= S^{-1}(\frac{U+S(\xu)}{2})\subset \tilde {\cal O}$; let $V_{\C}$ denote the subset of $\C$ given by the same formula. If $\Im S(\xu) < \Im S(x_1)$, fig. \ref{Paper3p44n},b),  then it is possible to draw an integration path $y(t)$  in $\tilde{\cal O}$ from $x_0$ and $\xu$. In a careful treatment, one sees that the set $U$ can indeed be transported along $y(t)$ parallel to $-S_1(x)$. 

If, on the contrary, $\Im S(\xu) > \Im S(x_1)$ and $N$ is large enough, $x_0$ and $\xu$ belong to different connected component of $\tilde {\cal O}\backslash V_\C$, fig. \ref{Paper3p44n},c), and an integration path $y(t)$ cannot be drawn. The situation is however remedied if instead of $U$ one considers a smaller subset $U'=U\backslash B_\varepsilon(2S(x_1)-S(\xu)+\R_{\ge 0})$ where $B_\varepsilon$ denotes an $\varepsilon$-neighborhood of a subset of the complex plane of $s$, for $\varepsilon>0$ small enough; the set $U'$ and the corresponding set $V'=S^{-1}(\frac{U'+S(\xu)}{2})\subset \tilde{\cal O}$ and a possible path $y(t)$ are shown on fig.\ref{Paper3p44n},d). This strongly suggests that for $\Im S(x)> \Im S(x_1)$ the first sheet of the Riemann surface of $R_1G_2$ contains a singularity at $s=2S(x_1)-S(\xu)$. We immediately recognize the curve $\Im S(x)=\Im S(x_1)$ as the Stokes curve and appearance of the new singularity as the Stokes phenomenon known, e.g., from ~\cite{V83}.


\section{Riemann surface ${\cal S}$} \label{Notation}

The notation presented here is heavily inspired by the collaboration on ~\cite{GT}.

\subsection{Assumptions on $V(x)$.} \label{WhoIsVx}
 We assume $V(x)$ to be a nonzero entire function of one complex variable $x$. Let $x_1,...,x_j,... $ denote the zeros of $V(x)$; zeros of $V(x)$ are also called {\it turning points} of the equation \eqref{SchroeEq}.

Let $x_0\in \C$ and $V(x_0)\ne 0$. Let $\tilde {\cal O}$ be the universal cover of $\C\backslash V^{-1}(0)$ with the base point $x_0$ and let $P:\tilde{\cal O}\to \C$ be the projection; the distinguished preimage of $x_0$ in $\tilde {\cal O}$ will also be denoted by $x_0$. On $\tilde {\cal O}$, let us choose and fix a determination of $p(x)=\sqrt{V(x)}$; then the integral $S(x)=\int_{x_0}^x \sqrt{V(y)}dy$ defines an analytic function on $\tilde {\cal O}$. It will be convenient to write $p_1(x):=p(x)$, $p_2(x):=-p(x)$, $S_1(x):=S(x)$, $S_2(x):=-S(x)$. 

The function $S(x)$ descends as a multivalued function to the complex plane of $x$, but it does not prevent us from making the following definition.

\begin{Definition} {\rm A {\it Stokes curve on $\C$} starting at a turning point $x_j$ is a connected component of the real curve defined by the condition $\Im [S(x)-S(x_j)]=0$, $x\ne x_j$. \\
A real curve $\ell$ on $\tilde {\cal O}$ is said to be a Stokes curve if its projection $P(\ell)$ to $\C$ is a Stokes curve. 
}
\end{Definition}

A Stokes curve on $\C$ originates at a turning point, but there are no preimages of turning points in $\tilde {\cal O}$. Nevertheless, we will, by abuse of language, say that a Stokes curve $\ell$ on $\tilde {\cal O}$ starts at a turning point $\tilde x_j$ in $\tilde {\cal O}$. Since, in fact, we will only be using the value $S(\tilde x_j)$ which can be unambiguously defined as a limit, this way of speaking will not lead to logical mistakes.

The following two assumptions on $V(x)$ will help us make the expositon shorter:

\begin{Assumption} \label{nonzeroLC} {\rm There is no nonzero collection of integers $m_1,...,m_j,...$, with only finitely many of them nonzero, such that 
$$ \Im [ m_1 S(x_1) + ... + m_j S(x_j) + ... ] \ = \ 0. $$
}\end{Assumption}

\begin{Assumption} \label{unbdd} {\rm All Stokes curves are unbounded; in other words, no Stokes curve connects two turning points.
}\end{Assumption}

The assumptions on $V(x)$ have now been listed completely.

\subsection{Stokes regions}

The Stokes curves in $\tilde{\cal O}$ split $\tilde{\cal O}$ into open regions which we will call the {\it Stokes regions}.  Let ${\mathbb S}$ denote the set of all Stokes regions. Let ${\cal A}_0$ be the Stokes region containing $x_0$. 

There is a partial order on ${\mathbb S}$: we say that a Stokes region ${\cal B}$ is {\it closer to $x_0$} than ${\cal B}'$ if every curve connecting $x_0$ to a point in ${\cal B}'$ has to pass through ${\cal B}$. Clearly, ${\cal A}_0$ is closer to $x_0$ than any other Stokes region. 

\subsection{Fiber of ${\cal S}$ over $x$} \label{FibOfS}

In this section we will describe ${\cal S}$ as a subset of $\tilde {\cal O}\times \C$ with the induced structure of a complex two-dimensional manifold. For each $x\in \tilde {\cal O}$, we will specify which horizontal rays $c+\R_{\ge 0}$ should be removed from $\C$ in order to obtain the fiber ${\cal S}_x$ of ${\cal S}$ over $x$. 

 Let ${\cal L}_S$ denote the set of all Stokes curves lifted to the universal cover. Let ${\cal L}_{\pravo}$ be the set of those Stokes curves for which $\Re S$ grows away from the turning point, ${\cal L}_{\levo}$ the set of those for which $\Re S$ grows towards the turning point, so that ${\cal L}_S = {\cal L}_\pravo \cup {\cal L}_\levo$.

The definition of ${\cal S}_x$ will proceed by induction. For ${\cal B}\in {\mathbb S}$, denote by $\bar{\cal B}$ the closure of ${\cal B}$ in $\tilde {\cal O}$. We will define ${\cal S}_x$ when $x\in \bar{\cal B}$ separately for each ${\cal B}\in {\mathbb S}$. Our definitions will agree for $x$ on the Stokes curves which are common boundaries of two Stokes region. 

If $x\in \bar{\cal B}$, we will define two sets of functions $\Sigma^+_{\cal B}$ and $\Sigma^-_{\cal B}$, $\bar{\cal B}\to \C$, and put $\Sigma_{\cal B}=\Sigma^+_{\cal B}\cup \Sigma^-_{\cal B}$ and let ${\cal S}_x:= \C\backslash \bigcup_{\sigma\in \Sigma_{\cal B}} (\sigma(x)+\R_{\ge 0})$. 

We will refer to $\sigma(x)\in \Sigma_{\cal B}$ as {\it singularities on the first sheet of ${\cal S}$ over the Stokes region ${\cal B}$.} Intuitively, if one thinks of an endless ``Riemann surface" of a solution of \eqref{MainEqu1} , then $s=\sigma(x)$ should represent the ramification curves of this Riemann surface.

In our definition, the functions in $\Sigma^\pm_{\cal B}$ will be of the form $const\pm S(x)$.

Set $\Sigma^+_{{\cal A}_0} := \{ S(x) \}$ and set $\Sigma^{-}_{{\cal A}_0}:=\{ -S(x)\}$.

Suppose we have defined the sets $\Sigma^{\pm}_{\cal B}$; let us define $\Sigma^{\pm}_{\cal C}$ where ${\cal C}$ is the Stokes region farther away from $x_0$ than ${\cal B}$ and such that $\bar{\cal B}\cap \bar{\cal C}=\ell\in {\cal L}_S$. Suppose $\ell$ ``starts at the turning point $x_t$ on $\tilde {\cal O}$", i.e. we will use the limit values of $\sigma(x)$ for $\sigma\in \Sigma_{\cal B}$ as $x$ tends to the origin of $\ell$ and denote them by $\sigma(x_t)$.

If $\ell\in {\cal L}_\levo$, then put $\Sigma^+_{\cal C}:=\Sigma^+_{\cal B}$ and $\Sigma^-_{\cal C}:=\Sigma^-_{\cal B} \cup \{ 2\sigma(x_t)-\sigma(x) : \sigma\in \Sigma^+_{\cal B}\}$. 

If $\ell\in {\cal L}_\pravo$, then put $\Sigma^-_{\cal C}:=\Sigma^-_{\cal B}$ and $\Sigma^+_{\cal C}:=\Sigma^+_{\cal B} \cup \{ 2\sigma(x_t)-\sigma(x) : \sigma\in \Sigma^-_{\cal B}\}$. 

These definitions are a reformulation of the description of the Stokes phenomenon in ~\cite{V83}.

Finally, let $\pi:{\cal S}\to \tilde{\cal O}$ denote the obvious projection. \label{piDefd}

The definition of ${\cal S}$ in now complete.

\subsection{Virtual Stokes curves and Stokes subregions}

Let ${\cal B}$ be a Stokes region. A point $x_v\in \bar{\cal B}$ is called a {\it virtual turning point} if there are singularities $\sigma\in \Sigma^+_{\cal B}$ and $\sigma'\in \Sigma^-_{\cal B}$ such that $\sigma(x_v)=\sigma'(x_v)$. 

It follows from assumption \ref{nonzeroLC} that all virtual turning poins are contained in the interior of the Stokes regions.

The curve in ${\cal B}$ defined by the equation $\Im [S(x)-S(x_t)]=0$ and unbounded in both directions is called a {\it virtual Stokes curve}. There are only finitely many virtual Stokes curves in every Stokes region.

The virtual Stokes curves split the Stokes region into open {\it Stokes subregions}. Let $Sub{\mathbb S}$ denote the set of all Stokes subregions in ${\cal S}$. Stokes subregions will be partially ordered with respect to their distance to $x_0$.

We will freely use the notation $\Sigma^{\pm}_{\cal B}$ when ${\cal B}$ is a Stokes subregion. In this case, if $\sigma_1,\sigma_2\in \Sigma_{\cal B}$ and $\Im \sigma_1(x) < \Im \sigma_2(x)$ for one point $x\in {\cal B}$, then the same inequality holds true for all points in ${\cal B}$.

The following lemma is immediate from our definitions. It is important that $\ell$ in this lemma is an actual, not a virtual Stokes curve.

\begin{Lemma} \label{noConsPairs} For ${\cal B}$ and ${\cal C}$ be Stokes subregions separated by a Stokes curve $\ell$ and ${\cal B}$ is closer to $x_0$ than ${\cal C}$. \\
i) If $\ell\in {\cal L}_\levo$, $\sigma_1,\sigma_3\in \Sigma^+_{\cal C}$ with $\Im \sigma_1<\Im \sigma_3$ in ${\cal C}$, then there is $\sigma_2\in \Sigma^-_{\cal C}$ with $\Im \sigma_1 < \Im \sigma_2 < \Im \sigma_3$. \\ 
ii) If $\ell\in {\cal L}_\pravo$, $\sigma_1,\sigma_3\in \Sigma^-_{\cal C}$ with $\Im \sigma_1<\Im \sigma_3$ in ${\cal C}$, then there is $\sigma_2\in \Sigma^+_{\cal C}$ with $\Im \sigma_1 < \Im \sigma_2 < \Im \sigma_3$.
\end{Lemma}

\subsection{Strips and slots.} \label{FSS}

Let $k\in \Z_{\ge 0}$, $s_1,...,s_k \in \C$ and $\Im s_{j+1} < \Im s_{j}$. 
 Consider the set $U_0$ obtained from $\C$ by removing horizontal cuts starting at $s_1,...,s_k$:
$$ U_0 = \C \backslash \bigcup_{j=1}^k (s_k + \R_{\ge 0}). $$
Every ${\cal S}_x$ is of the form $U_0$ for a suitable choice of $s_1,...,s_k$.
 
\begin{figure} \includegraphics{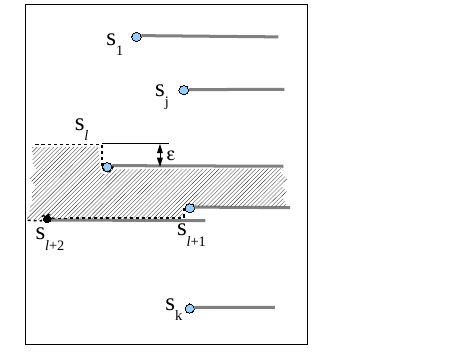} \caption{The $\varepsilon$-strip between $s_\ell$ and $s_{\ell+1}$} \label{Paper3p48} \end{figure}

\paragraph*{Strips.}
For a number $\varepsilon>0$, we are going to define a subset of $U_0$ which we will call the {\it $\varepsilon$-strip between $s_\ell$ and $s_{\ell+1}$}, fig.\ref{Paper3p48}, $1\le \ell \le k-1$ as the intersection of $A_{\ell+1}\cap B_{\ell}$,
$$ \begin{array}{clc} A_\ell=\{ s\in \C \ : & \Im s > \Im s_\ell - \varepsilon \\ 
                                            & (\Re s\ge \Re s_\ell) \Rightarrow (\Im s > \Im s_{\ell}) \\
                                            & ( \Re s< \Re s_\ell \ \text{and} \ \Re s_{\ell'}<\Re s_{\ell} \ \text{and} \ \ell'>\ell ) \Rightarrow (\Im s>\Im s_{\ell'}) & \} \end{array}; $$
$$ \begin{array}{clc} B_\ell=\{ s\in \C \ : & \Im s < \Im s_\ell + \varepsilon \\ 
                                            & (\Re \ge \Re s_\ell) \Rightarrow (\Im s < \Im s_{\ell}) \\
                                            & ( \Re s< \Re s_\ell \ \text{and} \ \Re s_{\ell'}<\Re s_{\ell} \ \text{and} \ \ell'<\ell ) \Rightarrow (\Im s<\Im s_{\ell'}) & \} \end{array}. $$
The {\it (semi-infinite) $\varepsilon$-strip above $s_1$} is defined to be $A_1$, the {\it (semi-infinite) $\varepsilon$-strip below $s_k$} is defined to be $B_k$. 


The prefix ``$\varepsilon$-" will sometimes be omitted.

\paragraph*{Slots.} For $s_0\in\C$, $\rho\ge 0$, define the {\it upward-facing slot of size $\rho$ around $s_0$}
$$Sl^\cup_\rho = \{ s\in \C : |s-s_0|\le \rho \} \cup \{ s\in \C : |\Re (s - s_0)|\le \rho \ \text{and} \ \Im s \ge \Im s_0 \} $$ 
and the {\it downward-facing slot of size $\rho$ around $s_0$} as
$$Sl^\cap_\rho = \{ s\in \C : |s-s_0|\le \rho \} \cup \{ s\in \C : |\Re (s - s_0)|\le \rho \ \text{and} \ \Im s \le \Im s_0 \}. $$



\section{Construction of analytic continuations.} \label{PfMresSec}


\subsection{Strategy of the proof. } \label{StrategySec}

The definitions  in the section  \ref{Notation}  have given a precise sense to content of the section \ref{RedToPaths}; we are continuing now where we stopped at the end of the section \ref{RedToPaths}. 

In order to prove theorem  \ref{RjGmainTh}, we will construct the integration paths for $R_j$ and $(s,\xu)\in {\cal S}$ from $x_0$ to $\xu$ by induction on the Stokes subregions containing $\xu$ in their closure. We will first consider the case when $\xu$ belongs to one of the two Stokes subregions ${\cal A}'_0$, ${\cal A}''_0$ which contain $x_0$ on their boundary. In order to proceed by induction with respect to the partial order on the set $Sub{\mathbb S}$ of all Stokes subregions, for $(s,\xu)\in{\cal S}_\xu$ with $\xu\in {\cal C}\in Sub{\mathbb S}$ we will construct a piece of integration path that starts at $\xu'\in {\cal B}\in Sub{\mathbb S}$ and leads to $\xu$, with ${\cal B}$ closer to $x_0$ than ${\cal C}$. 
\label{IndArgProposed}

\subsection{Base of induction}
\label{ProofAint}


Recall that ${\cal A}_0$ denotes the Stokes region containing the point $x_0$. Consider the virtual Stokes curve $L_0$ in ${\cal A}_0$ given by the equation $\Im S(x)=0$ passing through the point $x_0$; it splits ${\cal A}_0$ into two Stokes subregions ${\cal A}'_0$ where $\Im S(x)>0$ and ${\cal A}''_0$ where $\Im S(x)<0$. 

In this section \ref{ProofAint}, we will consider ${\cal A}'_0$ and discuss a construction of $R_1G$; the other three pairs of choices between ${\cal A}'_0$ and ${\cal A}''_0$ and between $R_1$ and $R_2$ are analogous.

\begin{figure} \includegraphics{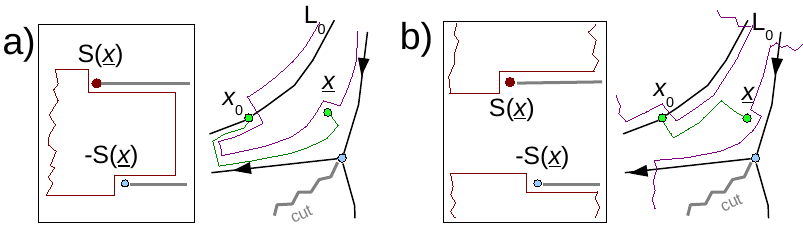} \caption{ The sets  a) $U=U^{\eta}_{\xu,N}$, b) $U=W^\eta_\xu$ in ${\cal S}_\xu$;  the corresponding sets $V$ in $\tilde{\cal O}$, and the integration paths $y(t)$ from $x_0$ to $\xu$. The arrows on the Stokes curves show the direction in which $\Re S$ grows.} \label{Paper3p55n} \end{figure}

Let $\xu\in {\cal A}'_0\cup L_0$. Let us write ${\cal S}_\xu$ as a union of more convenient sets. Fix some $\varepsilon>0$. For every (large) $N\in \R$, let $U_{\xu,N}$ be the intersection of the $\varepsilon$-strip between $S(\xu)$ and $-S(\xu)$ with the set $\{s :\, \Im s<N\}$. Let $W_{\xu}$ be the union of the semi-infinite $\varepsilon$-strip above $S(\xu)$ and the semi-infinite $\varepsilon$-strip below $-S(\xu)$.  For every (small) $\eta>0$, let 
$$U^\eta_{\xu,N} = ( U_{\xu,N} - i\eta) \cap ( U_{\xu,N} + i\eta) \cap (U_{\xu, N}- \eta), $$
$$W^\eta_\xu = ( W_\xu - i\eta) \cap ( W_\xu + i\eta) \cap ( W_\xu - \eta). $$
For each fixed $\varepsilon>0$, we have
$$ {\cal S}_x = \bigcup_{N,\eta>0} ( U^\eta_{\xu,N} \cup W^\eta). $$
Let $U$ be one of the sets $U^\eta_{\xu,N}$ or $W^\eta_x$. We will now present the integraton path $y(t)$ for the points $(s,\xu)$ where $s\in U$. Consider $V=S^{-1}\left( \frac{U+S(x)}{2} \right)$ where the branch of $S^{-1}$ is chosen in such a way that $S^{-1}(x)=x_0$, fig. \ref{Paper3p55n}. We have made our definition in such a way that the set $({\cal A}'_0 \cup L_0)\backslash V$ is connected, and there is a path $y(t)$ from $x_0$ to $\xu$ in $({\cal A}'_0 \cup L_0)\backslash V$ which can be taken as an integration path, fig. \ref{Paper3p55n}.

\label{WhyAnal}
It is easy to see that this defines an {\it analytic} function in a neighborhood of a point $(s,\xu)$, with $\xu\in{\cal A}'_0$, $s\in U$. Indeed, if $\xu'$ is another point in ${\cal A}'_0$, then the set $U$ can also be transported parallel to $-S_1$ along any path $z(t)$ from $\xu$ to $\xu'$ if $z(t)$ is contained in ${\cal A}'_0\backslash V$. There is an open contractible set ${\cal N}_\xu\subset \tilde {\cal O}$ such that $\xu\in {\cal N}_\xu\subset {\cal A}'_0\backslash V$.  For any $s'\in U$ and $x\in {\cal N}_\xu$, the function  $(R_1G)(s'+S_1(\xu)-S_1(x),x)$ is holomorphic with respect to $x$ because it is an integral of a holomorphic function and with respect to $s'$ because $s'$ is a holomorphic parameter of the integrand, hence, by Osgood theorem, this function is holomorphic in both $s'$ and $x$, and so, after a change of variables, is $(R_1G)(s',x)$ is holomorphic in $(s',x)$ in a neighborhood of $(s,\xu)$.

The reader will easily make an argument along the same lines for $\xu \in L_0$.

Recall, p.\pageref{piDefd}, that $\pi:{\cal S}\to\tilde {\cal O}$ denotes the obvious projection.

We have shown: \\
{\it If $G(s,x)$ is an analytic function on $\pi^{-1}{\cal A}_0\subset {\cal S}$, then so is $(R_jG)(s,x)$, $j=1,2$. }


\subsection{Continuation to the further Stokes subregions.} \label{modelcases}

The rest of the argument will procede by induction on the set $Sub{\mathbb S}$ -- we will construct $R_1G(s,x)$ for $x$ in the Stokes subregions further and further away from $x_0$. 

Let ${\cal B},{\cal C}$ be two Stokes subregions, ${\cal B}$ closer to $x_0$ than ${\cal C}$ and $\bar{\cal B}\cap \bar{\cal C}=\ell$, where $\ell$ is either a Stokes curve or a virtual Stokes curve.

For either choice of $j=1,2$ and for every point $(s,x)\in {\cal S}$, $x\in {\cal C}\cup \ell$, we will present an piece of an integration path from $x$ to some point in ${\cal B}$. This will accomplish the analytic continuation of $R_jG$ to the Stokes region ${\cal C}$.

Here and in the next section \ref{StokesXing} we will discuss in detail the case when $\ell$ is a Stokes curve in $\tilde {\cal O}$  ``starting at the turning point $x_t$"; in the section \ref{VirtualXing} we will make a short remark what changes when $\ell$ is a virtual Stokes curve.

We will assume that the order of ${\cal B},\ell,{\cal C}$ around $x_t$ is clockwise; if the order is clockwise, one should exchange the positive imaginary and the negative imaginary directions in the $s$-plane in all the statements below.


Recall that when constructing an integration path for $R_j$, we call singularities of the form $-S_j(x)+const$ {\it stationary} and those of the form $S_j(x)+const$ {\it moving}.

The lemmas we are going to present now can be systematized in the following manner. Fistly, each lemma may pertain to a case when $\Re S_j$ increases, resp., decreases along $\ell$. Secondly, we will choose our set $U$ (notation of p.\pageref{Ucomes}) as a strip between two singularities: moving or stationary above and moving or stationary below ($2\times 2=4$ possibilities), or as a semi-infinite strip, above or below a moving or a stationary singularity ($2\times 2=4$ more possibilities). Multiplying this by $2$ to account for the direction of growth of $\Re S_j$, we obtain $(4+4)\times 2=16$ possibilities. The reader will see that this is an easily manageable number of cases.

\subsection{Crossing actual Stokes curves} \label{StokesXing}

We will use the following Assumption in every lemma of this subsection \ref{StokesXing}.

\begin{Assumption} \label{XlY} Let ${\cal B},{\cal C}$ be two Stokes subregions, ${\cal B}$ closer to $x_0$ than ${\cal C}$, $\bar{\cal B}\cap\bar{\cal C}=\ell$ be a Stokes curve.  Let $G(s,x)$ be an analytic function on ${\cal S}$.
\end{Assumption}

We will say that $S_j$ {\it grows, resp., decays along the Stokes curve $\ell$} if $\Re S_j(x)$ increases, resp., decreases as $x$ moves away from the turning point.

If ${\cal C}$ is a Stokes subregion, we will say that $\sigma_1,...,\sigma_m \in \Sigma_{\cal C}$ are {\it consecutive} singularities in ${\cal C}$ if for every $k=1,..,m-1$ and for some (hence any) point $x\in {\cal C}$ we have $\Im \sigma_k(x)<\Im \sigma_{k+1}(x)$ and there is no $\sigma'\in \Sigma_{\cal C}$ with $\Im \sigma_k(x) <\Im \sigma'(x) < \Im \sigma_{k+1}(x)$.

As before, $\pi:{\cal S}\to \tilde{\cal O}$ denotes the obvious projection. We will often and freely identify a subset $D\subset {\cal S}$ with a subset of $\C\times \tilde {\cal O}$ or with a subset of $\C\times \C$; we hope that this will not cause any confusion. We will use an abbreviation ``a function $G$ is C.A.I. in $D$" to mean that $G$ is continuous on a set $D$ and analytic in its interior.

With this, let us start working through the 16 cases mentioned at the end of section \ref{modelcases}.

\subsubsection{A strip between two moving singularities, $S_j$ grows along $\ell$} 

This situation is never realized, see Lemma \ref{noConsPairs}.

\subsubsection{A strip between a stationary singularity (above) and a moving singularity (below), $S_j$ grows along $\ell$} \label{SeOc4}

Under assumptions \ref{XlY}, given consecutive singularities $\sigma_2,\sigma_3\in \Sigma_{\cal C}$, $\sigma_2$ moving, $\sigma_3$  stationary, there must be, by the inductive definition of $\Sigma_{\cal C}$ in terms of $\Sigma_{\cal B}$, a stationary singularity $\sigma_1=2(\sigma_2(x_t)-\sigma_2) \in \Sigma_{\cal C}$. This puts us in the situation of the following lemma.

\begin{Lemma} \label{YelGrIN} Under assumptions \ref{XlY}, let $S_j$ grow along $\ell$, let $\sigma_1,\sigma_2,\sigma_3\in \Sigma_{\cal C}$ be consecutive singularities in ${\cal C}$, where  $\sigma_1$, $\sigma_3$ are stationary, $\sigma_2$ is moving, $\sigma_1(x_t)=\sigma_2(x_t)$. Consider a function $\sigma_4(x)=2\sigma_3(x_t)-\sigma_3(x)$ which may or may not belong to $\Sigma_{\cal C}$. Fix an $\varepsilon>0$. 
\\
Suppose the function $R_jG(s,x)$ is C.A.I. in the set
$$ \{ (s,x)\in {\cal S} \ : \ x\in {\cal B}, \text{ $s$ in the $\varepsilon$-strip between $\sigma_1$ and $\sigma_4$} \ \}. $$
Then $R_jG$ has an analytic continuation to
$$ \{ (s,x)\in {\cal S}\ : \ x\in \bar{\cal C}, \text{$s$ in the $\varepsilon$-strip between $\sigma_2$ and $\sigma_3$}\}. $$
\end{Lemma}

\textsc{Proof.} Let us note first that in $\Sigma_{\cal B}$ there may (or may not) be a moving singularity $\sigma_4$ such that $\sigma_4(x_t)=\sigma_3(x_t)$. Also, note that $\sigma_1$ may or may not be present in the set $\Sigma_{\cal B}$.  We will carry out the proof in the case when both $\sigma_1,\sigma_4\in \Sigma_{\cal B}$; the other cases are similar but simpler.

Under the above assumption,  $\sigma_4,\sigma_3,\sigma_1$ are consecutive singularities in $\Sigma_{\cal B}$.

\begin{figure} \includegraphics{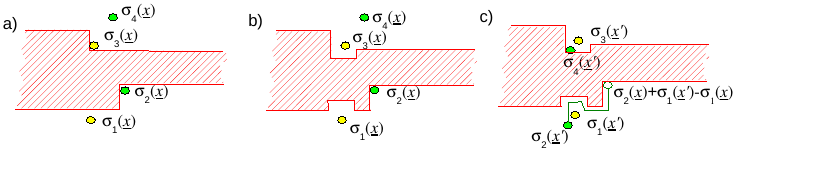} \caption{Proof of lemma \ref{YelGrIN}. a) The set $B^{\xu}_{\sigma_2\sigma_3}\subset {\cal S}_\xu$. b) The set $B'\subset {\cal S}_{\xu}$ c) The set $B'+\sigma_1(\xu')-\sigma_1(\xu) \subset {\cal S}_{\xu'}$ and the path $\tau$. } \label{Paper3p51} \end{figure}

 Denote, for $\xu\in \bar {\cal C}$, by $B^{\xu}_{\sigma_2\sigma_3}$ the $\varepsilon$-strip between $\sigma_2(\xu)$ and $\sigma_3(\xu)$, see fig.\ref{Paper3p51},a). This set $B^{\xu}_{\sigma_2\sigma_3}$ can be exausted by its subsets  
$$ B' = B^\xu_{\sigma_2 \sigma_3} \backslash [ Sl^\cup_\eta(\sigma_3(\xu)) \cup Sl^\cap_\eta(\sigma_1(\xu)) ].$$
for parametrized by $\eta$, $0<\eta<\varepsilon/2$, see fig.\ref{Paper3p51},b). It is thus enough to construct the analytic continuation of $R_1G$ to the points $(s,\xu)$ for $s\in B'$ for each fixed $\eta$.

Take any $\xu'$ in ${\cal B}$ such that $\sigma_4(\xu')\in Sl^\cup_\eta(\sigma_3(\xu'))$. Consider a path $\tau$ in the $\C\backslash [B'+\sigma_1(\xu')-\sigma_1(\xu)]$ from $\sigma_2(\xu')$ to $\sigma_2(\xu)+\sigma_1(\xu')-\sigma_1(\xu)$, going clockwise around $\sigma_1(\xu')$, fig.\ref{Paper3p51},c) 

Consider the preimage of this path $\sigma_2^{-1}(\tau)$ in ${\cal B}\cup \ell \cup {\cal C}$. The set $B'$ can be  transported parallel to $-S_j$ (cf. Def.\ref{cbpted}) along the path $\sigma_2^{-1}(\tau)$ and the function, by assumptions of the lemma,  $R_1G(s,x)$ is defined for $x=\xu'$, $s\in B'+\sigma_1(\xu')-\sigma_1(\xu)$. Thus we can use $\sigma_2^{-1}(\tau)$ as an integration path to obtain the value of $R_1G(s-\sigma_1(\xu')+\sigma_1(\xu),\xu)$.

This definition clearly gives an analytic continuation of $R_jG$. $\Box$

\subsubsection{A strip between two stationary singularities, $S_j$ grows along $\ell$} 

This case is similar to section \ref{SeOc4}. We will let the reader formulate the corresponding lemma. 

\subsubsection{A strip between a moving singularity (above) and a stationary singularity (below), $S_j$ grows along $\ell$} 

Suppose we have constructed the analytic continuation of $R_1G(s,\xu)$ for $\xu$ in the Stokes subregion ${\cal C}$ for all $\varepsilon$-strips in ${\cal S}_\xu$ between a stationary singularity above the strip and a moving singularity below the strip, as in section \ref{SeOc4}. 

Let us now consider the strip in ${\cal S}_\xu$ between two consecutive singularities $\sigma_1,\sigma_2\in \Sigma_{\cal C}$, where $\sigma_1$ is stationary and $\sigma_2$ is moving. Suppose there is yet another singularity $\sigma_3$ such that $\sigma_1,\sigma_2,\sigma_3\in \Sigma_{\cal C}$ are consecutive; if there is no such $\sigma_3$, the argument below will work once we replace $\sigma_3$ with $+i\infty$. 

By lemma \ref{noConsPairs}, $\sigma_3$ is a stationary singularity. Thus, and in view of section \ref{SeOc4}, we are in the situation of the following lemma. 

%

\begin{Lemma} \label{GrYelIN} Under assumptions \ref{XlY}, let $S_j$ grow along $\ell$. Let $G(s,x)$ be defined on ${\cal S}$. Let $\sigma_1, \sigma_2, \sigma_3 \in \Sigma_{\cal C}$ be consecutive singularities, $\sigma_1,\sigma_3$ stationary, $\sigma_2$ moving, $\sigma_1(x_t)=\sigma_2(x_t)$.  Let $\varepsilon>0$.
Suppose the function $R_j G(s,x)$ is C.A.I. in the set 
$$ \{ (s,x)\in {\cal S} \ : \ x\in {\cal C}, \ s \, \text{in the $\varepsilon$-strip between $\sigma_2$ and $\sigma_3$} \}. $$
Then $R_j G(s,x)$ has an analytic continuation to 
$$ \{ (s,x)\in {\cal S} \ : \ x\in \bar{\cal C}, \ s \, \text{in the $\varepsilon$-strip between $\sigma_1$ and $\sigma_2$} \}. $$
\end{Lemma} 

\begin{figure} \includegraphics{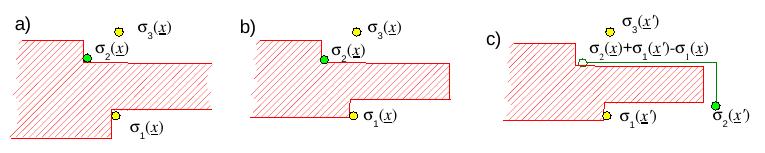} \caption{Proof of Lemma \ref{GrYelIN}. a) The set $B^{\xu}_{\sigma_1\sigma_2}\subset {\cal S}_\xu $; b) The set $B'\subset {\cal S}_\xu$; c) the set $B'+\sigma_1(\xu')-\sigma_1(\xu)\subset {\cal S}_{\xu'}$ and the path $\tau$.} \label{Paper3p50} \end{figure}

\textsc{Proof.}  Denote, for $x$ in the closure of ${\cal C}$, by $B^{x}_{\sigma_1,\sigma_2}$ the $\varepsilon$-strip between $\sigma_1(x)$ and $\sigma_2(x)$, see fig.\ref{Paper3p50},a), and similarly for $B^x_{\sigma_2,\sigma_3}$. 

Let $\xu$ be in the closure of ${\cal C}$. 
Let $N>0$, let $\eta>0$. We can exhaust $B^\xu_{\sigma_1 \sigma_2}$ by the subsets of the form 
$$ B'=B^x_{\sigma_1 \sigma_2} \cap (B^x_{\sigma_1 \sigma_2}+i\eta) \cap \{ \Re s < \Re \sigma_1(x)+N\}, $$
thus it is enough to construct the analytic continuation of $R_1G(s,x)$ to every set $B'$ for every sufficiently large $N$ and every sufficiently small $\eta$.

Let us find $\xu'\in Y$ such that $\sigma_2(\xu')-\sigma_1(\xu') = N+1 + i\frac{\eta}{2}$. 

Consider a two line segment path $\tau$ in the $s$-plane starting at $\sigma_2(\xu')$, ending at $\sigma_2(\xu)+\sigma_1(\xu')-\sigma_1(\xu)$ as on fig. \ref{Paper3p50},c).  The set $B'$ can be parallel transported along the path $\sigma_2^{-1}(\tau)$ in $\bar{\cal C}$  parallel to $-S_j$ from $\xu$ to $\xu'$, and the function $R_jG(s+\sigma_1(x)-\sigma_1(\xu),x)$ is defined for $x=\xu'$ and $s\in B'$ because $B'+\sigma_1(\xu')-\sigma_1(\xu)\subset B^{\xu'}_{\sigma_2 \sigma_3}$ by assumptions of the lemma. Thus we can use  $\sigma_2^{-1}(\tau)$ as an integration path to obtain the value of $R_jG(s,\xu)$. This definition clearly gives an analytic continuation of $R_j G$. $\Box$


\subsubsection{Semiinfinite strips, $S_j$ grows along $\ell$}

The cases of \\
{\it i)} a semi-infinite strip between $+i\infty$ and a stationary singularity, \\
{\it ii)} a semi-infinite strip between $-i\infty$ and a stationary singularity, \\
{\it iii)} a semi-infinite strip between $+i\infty$ and a moving singularity, \\
are treated similarly to lemma \ref{GrGrOUT}. 

The case of a semi-infinite strip between $-i\infty$ and a moving singularity is never realized, compare Lemma \ref{noConsPairs}.

\subsubsection{A strip between two moving singularities, $S_j$ decays along $\ell$}

\begin{Lemma} \label{GrGrOUT}  Under assumptions \ref{XlY}, let $S_j$ decay along $\ell$. Let $\sigma_2,\sigma_3\in \Sigma_{\cal C}$ be consecutive singularities, both of them moving. Let $\varepsilon>0$. \\
Suppose the function $R_jG(s,x)$ is defined and analytic for $x\in {\cal B}$. 
Then $R_jG$ has an analytic continuation to
$$ \{ (s,x)\in {\cal S}\ : \ x\in \bar{\cal C}, \text{$s$ in the $\varepsilon$-strip between $\sigma_2$ and $\sigma_3$}\}. $$
\end{Lemma}

\begin{figure} \includegraphics{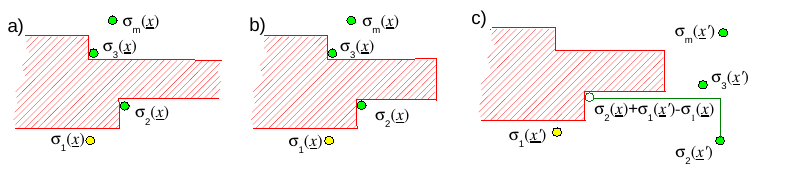} \caption{Proof of Lemma \ref{GrGrOUT}.  a) The set $B^{\xu}_{\sigma_2\sigma_3}\subset {\cal S}_\xu $; b) The set $B'\subset {\cal S}_\xu$; c) the set $B'+\sigma_2(\xu')-\sigma_2(\xu)\subset {\cal S}_{\xu'}$ and the path $\tau$.} \label{Paper3p52} \end{figure}

\textsc{Proof.}  Consider the function $\sigma_1(x)=2\sigma_2(x_t)-\sigma_2(x)$  which may or may not belong to $\Sigma_{\cal C}$. (On Fig.\ref{Paper3p52} we assume that $\sigma_1\in\Sigma_{\cal C}$.)

Let $B^{\xu}_{\sigma_2\sigma_3}$ denote the $\varepsilon$-strip between $\sigma_2(\xu)$ and $\sigma_3(\xu)$ in ${\cal S}_\xu$, fig.\ref{Paper3p52},a). Exhaust $B^{\xu}_{\sigma_2\sigma_3}$ by its subsets $B'=B^{\xu}_{\sigma_2\sigma_3}\cap \{ s: \Re s <  Re \sigma_1(\xu) + N\}$ for numbers $N>0$, fig.\ref{Paper3p52},b). Having fixed such a $B'$, choose $\xu'\in {\cal B}$ such that $\Re [\sigma_m(\xu')-\sigma_1(\xu)]>N$ for all moving singularities $\sigma_m\in \Sigma_{\cal C}$, and choose a path $\tau$ as a two segment broken line connecting $\sigma_2(\xu')$ and $\sigma_2(\xu)+\sigma_1(\xu')-\sigma_1(\xu)$, fig.\ref{Paper3p52},c). Take $\sigma_2^{-1}(\tau)$ as the integration path. This provides the desired analytic continuation of $R_jG$ by the same argument as in Lemma \ref{YelGrIN}. $\Box$

\subsubsection{A strip between a moving singularity (above) and a stationary singularity (below), $S_j$ decays along $\ell$}

The treatment of this case is very similar to Lemma \ref{GrGrOUT}. 

\subsubsection{A strip between a stationary singularity (above) and a moving singularity (below), $S_j$ decays along $\ell$}

This case is also similar to Lemma \ref{GrGrOUT}.

\subsubsection{A strip between two stationary singularities, $S_j$ decays along $\ell$}

This case is never realized by Lemma \ref{noConsPairs}.

\subsubsection{Semiinfinite strips, $S_j$ decays along $\ell$.}

The cases of \\
{\it i)} a semi-infinite strip between $+i\infty$ and a stationary singularity, \\
{\it ii)} a semi-infinite strip between $-i\infty$ and a stationary singularity, \\
{\it iii)} a semi-infinite strip between $-i\infty$ and a moving singularity, \\
are treated similarly to lemma \ref{GrGrOUT}. 

The case of a semi-infinite strip between $+i\infty$ and a moving singularity is never realized, compare Lemma \ref{noConsPairs}.   

\subsection{Crossing a virtual Stokes curve} \label{VirtualXing}

If the Stokes subregions ${\cal B}$ and ${\cal C}$ are separated by a virtual Stokes curve $\ell$ passing through a virtual turning point $x_v$, then the we can procede very analogously to the section \ref{StokesXing}. There are two main differences which turn out to compensate each other. On the one hand, we do not have a counterpart of Lemma \ref{noConsPairs} for a virtual Stokes curve. On the other hand, while actual Stokes curves are semi-infinite, virtual Stokes curves are infinite in both directions and an integration path can cross either of the two semi-infinite branches or even pass through the virtual turning point.

This finishes the proof of Theorem \ref{RjGmainTh}.

\vskip3pc

{\bf \large Acknowledgments} 

This work was mainly carried out during the author's studies at the Department of Mathematics, Northwestern University, U.S.A., and during his stay at the Max Planck Institute for Mathematics in the Sciences, Leipzig, Germany. The author is profoundly grateful to Dmitry Tamarkin and Boris Tsygan for numerous discussions and for conscientious critique of the manuscript, and also to Ovidiu Costin, Stavros Garoufalidis, Rostislav Matveyev, Shinji Sasaki, Boris Shapiro, Yoshitsugu Takei, and Jared Wunsch for discussions, comments, and feedback. The author also appreciates valuable comments of the anonymous referee.



\vspace{3cm}


\end{document}